\documentclass[10pt,leqno,twoside]{article}
\usepackage{amsmath,amssymb}
\author{ Rosa Cid-Mu$\tilde n$oz \and Manuel
Pedreira}
\title{Another Classification of Incidence Scrolls}
\date{}

\newtheorem{teo}{Theorem}[section]

\newtheorem{prop}[teo]{Proposition}

\newtheorem{lemma}[teo]{Lemma}

\newtheorem{rem}[teo]{Remark}

\newcommand{\p}{\ensuremath{\textrm{\rm I\hspace{-1.9pt}P}}}
\newcommand\om{\Omega}
\def\E{{\cal E}}
\def\Te{{\cal O}}

\font\euf=eufm10 at 14pt

\def\d{\mbox{\euf b}}

\def\b{{\cal B}}

\def\qed{\hspace{\fill}$\rule{1.5mm}{1.5mm}$}
\def\lrw{\longrightarrow}

\def\sub{\subset}

\begin{document}
\maketitle

{\footnotesize {\bf Authors' address:} Departamento de
Algebra, Universidad de Santiago de Compostela. $15782$
Santiago de Compostela. Galicia. Spain. Phone:
34-81563100-ext.13152. Fax: 34-81597054. {\tt e}-mail: {\tt
rosacid@usc.es}; {\tt pedreira@zmat.usc.es}\\ {\bf
Abstract:} The aim of this paper is the computation of the
degree and genus of  all incidence scrolls in $\p^n$. For
this, we fix the dimension of a linear space which must be
contained in the base of an incidence scroll. Then we will
find all the incidence scrolls which have a base space of
this fixed dimension. In this  way, we can obtain all the
incidence scrolls with a line as directrix curve, those
whose base contains a plane, and so on.}

\vspace{0.1cm}

{\bf Introduction:} Throughout this paper, the base field for algebraic varieties
is ${\mathbb{C}}$. Let $\p^n$ be the n-dimensional complex
projective space and $G(l,n)$ the Grassmannian of $l$-planes
in
$\p^n$. Then $R^d_g \sub \p^n$ denotes a scroll of degree $d$
and genus $g$. We will follow the notation and terminology of
\cite{hatshor}.

It is useful to represent a scroll in $\p^n$ by a curve
$C\sub G(1,n)\sub \p^N$. The lines which intersect
a given subspace $\p^r\sub \p^n$ are represented by
the points of the special Schubert variety $\om (\p^r,\p^n)$.
Each
$\om (\p^r,\p^n)$ is the intersection of $G(1,n)$ with a
certain subspace of $\p^N$. Since $G(1,n)$ has dimension
$2n-2$ and we search a curve, we must impose $2n-3$ linear  
conditions on $G(1,n)$. Consequently, the choice of subspaces
is not arbitrary. Any set of subspaces of $\p^n$ which
imposes $2n-3$ linear conditions on $G(1,n)$ is the base of
a certain incidence scroll. The background about Schubert
varieties can be found in \cite{Kleiman}.

The aim of this paper is to obtain another classification of
the scrolls in
$\p^n$ which are defined by a one-dimensional family of
lines meeting a certain set of subspaces of $\p^n$, a first
classification is given in paper
\cite{rosa}. These ruled surfaces  are called {\it
incidence scrolls}, and such an indicated set is a {\it
base} of the incidence scroll. Unless otherwise stated, we
can assume that the base spaces are in general position. We
will fix the dimension of a linear space which must be
contained in the base of an incidence scroll. In this  way,
we can obtain all the incidence scrolls which have a base
space of this fixed dimension, i.e., all the incidence
scrolls with a line as directrix curve, those whose base
contains a plane, and so on.

In section 1 we explain our notation and collect
background material. Our first step will be to summarize some
properties of ruled surfaces and, in par\-ticular, some
general properties of incidence scrolls. For more details we
refer the reader to \cite{hatshor} and \cite{rosa}
respectively. Having revised the notion of incidence scroll,
we have compiled some basic properties of such a scroll. This
section contains a brief summary and the detailed proofs
will appear in \cite{rosa}. The degree of the scroll given
by a general base is provided by Giambelli's formula which
appears in
\cite{gyh}. Moreover, the study of deformations of a given
incidence scroll is a powerful tool in this paper in order
to simplify our proofs. If the incidence scroll
$R^d_g
\sub \p^n$ breaks up into
$R^{d_{1}}_{g_{1}} \sub \p^r$ and $R^{d_{2}}_{g_{2}} \sub
\p^s$ with
$\delta$ generators in common, then $d=d_{1}+d_{2}$ and $
g=g_{1}+g_{2}+\delta-1$. 

Section 2 contains our main results about incidence scrolls.
We identify all the incidence scrolls. In this way, we 
find all incidence scrolls which contain a directrix line,
those which contain a base plane and finally those which
contain a base space of dimension 3. A more complete theory
may be obtained by a combination of these processes.
Accordingly, all the other incidence scrolls may be obtained.
Moreover, we will see that the fixed base space imposes
conditions on the genus of the incidence scroll.

The results on this paper 
belong to the Ph.D. thesis of the first author whose advisor 
is the second one.

\addtocontents{toc}{\protect\vspace{3ex}}
\section{Incidence Scrolls} \label{1}
\addtocontents{toc}{\protect\vspace{1ex}}

A {\it ruled surface} is a surface $X$ together with a
surjective morphism $\pi \colon X\lrw C$ to a
smooth curve $C$ such that the fibre $X_y$
is isomorphic to $\p^1$ for every point $y\in C$, and such
that $\pi$ admits a section. There exists a locally free
sheaf $\E$ of rank 2 on $C$ such that 
$X\cong \p (\E)$ over $C$. Conversely, every such $\p
(\E)$ is a ruled surface over $C$. 

A {\it scroll} is a ruled surface embedded in $\p^n$ in such
a way that the fibres $f$ have degree $1$. If we take a very
ample divisor on $X$, $D \sim aC_o+\d f$, then the
embedding $\Phi_D \colon X
\lrw R^d_g\sub\p^n=\p(H^0(\Te_X(D)))$ determines a scroll
when $a=1$. A scroll $R^d_g\sub
\p^n$ is said to be {\it an incidence scroll} if it is
generated by the lines which meet a certain set $\b$ of
linear spaces in $\p^n$, or equivalently, if the
correspondent curve in
$G(1,n)$ is an intersection of special Schubert varieties
$\om(\p^{r}, \p^{n})$, $0\leq r <n-1$. Such a set is called
a {\it base} of the incidence scroll and such a base will be
denoted by:
$$\b =\{\p^{n_1},\p^{n_2},\cdots ,\p^{n_r}\}.$$ 
We will write it simply $\b$ when no confusion can arise,
where $n_1\leq n_2\leq \cdots \leq n_r$, i.e.,
$C^d_g=\bigcap_{i=1}^{r}(\om(\p^{n_i}, \p^{n}))\sub G(1,n)$. 

Therefore,
unless otherwise stated, we will work with linear spaces in general position. By
general position we will mean that
$(\p^{n_1},\cdots ,\p^{n_r})\in {\cal W}=G(n_1,
n)\times \cdots \times G(n_r, n)$ is contained in
a nonempty open subset of
${\cal W}$. For simplicity of
notation, we abbreviate it to base in general
position. 
 
\begin{prop} \label{IS} The intersection
$C=\bigcap_{i=1}^{r}
\Omega(\p^{n_i},\p^{n})$ of special Schubert varieties
associated to linear spaces $\p^{n_i},\, i=1,\cdots,r$, in
general position is an irreducible curve of
$G(1,n)$ if and only if it verifies the following
equality
$$ rn-(n_1+n_2+\cdots +n_r)-r=2n-3 \eqno(IS) $$
\end{prop} 
{\bf Proof.} See \cite{rosa}, Proposition 2.4.\qed 

Moreover, the incidence scroll generated by a base $\b$
have degree $d$ if and only if we obtain the following equality of
Schubert cycles: $$\om(n_1, n) \cdots
\om(n_r, n)=d\, \om(0, 2).$$

Furthermore, we present one of the three main theorems of
the symbolic formalism, known as Schubert calculus, for
solving enumerative problems.

\begin{teo} \label{pieri}
(Pieri's formula) For all sequences of integers $0\leq
a_0<\cdots <a_l \leq n$ and for $h=0,\cdots,n-l$, the
following formula holds in the cohomology ring
$H^{\star}(G(l,n); {\mathbb{Z}})$:
$$
\om(a_0,\cdots,a_l)\om(h, n)=\sum \om(b_0,\cdots,b_l)
$$
where the sum ranges over all sequences of integers
$b_0<\cdots <b_l$ satisfying $0\leq b_0\leq a_0<b_1\leq
a_1<\cdots<b_l\leq a_l$ and
$\sum^l_{i=0}b_i=\sum^l_{i=0}a_i-(n-l-h)$.
\end{teo}
{\bf Proof.} See \cite{Kleiman}, p. 1073.\qed

Finally, let us mention an important property of
degeneration of these scrolls.

\begin{prop} \label{degen}
Let $R^d_g \sub \p^n$ be an incidence scroll with base 
$\b$ in general position. Suppose that $\p^{n_i}$ and
$\p^{n_j}$ lie in a  hyperplane $\p^{n-1}$ and have in
common
$\p^{m}$, 
$m=n_i+n_j-n+1$. Then
the scroll breaks up into: 
\begin{enumerate} 
\item[-] $R^{d_1}_{g{_1}} \sub \p^{n}$
with base $\dot \b
=\{\p^{m},\p^{n_1},\cdots,\widehat{\p^{n_i}},\cdots,\widehat{\p^{n_j}},
\cdots,\p^{n_r}\}$ (which is possibly degenerate);
\item[-] $R^{d_2}_{g_2} \sub \p^{n-1}$
with base $\ddot \b =\{\cdots ,
\p^{n_{i-1}-1},\p^{n_i},\p^{n_{i+1}-1},\cdots
,\p^{n_{j-1}-1},\p^{n_j},$ $\p^{n_{j+1}-1},\cdots \}$  
\end{enumerate}
which have $\kappa \geq 1$ generators in common. Then,
$d=d_1+d_2$ and $g=g_1+g_2+ \kappa -1$.

Moreover, if $m=0$, then the incidence scroll breaks up
into a plane and an incidence scroll $R^{d-1}_g \sub
\p^{n-1}$ with base $\ddot \b$ in general position. 
\end{prop}
{\bf Proof.} See \cite{rosa}, Proposition 3.1.\qed

If $m=0$, then shall refer to this particular degeneration as
{\it join
$\p^{n_i}$ and $\p^{n_j}$} (i.e., $n_i+n_j=n-1$) and to the
inverse as {\it separate
$\p^{n_i}$ and $\p^{n_j}$} (i.e., $n_i+n_j=n$).

\addtocontents{toc}{\protect\vspace{3ex}}
\section{Classification of Incidence Scrolls}
\label{2}
\addtocontents{toc}{\protect\vspace{1ex}}

In \cite{rosa} we have obtained a classification of
incidence scrolls of genus $0$ and $1$. We can study a new
point of view to give another classification of the incidence
scrolls. For this we will fix the dimension of a base space.
In this way, we will obtain all the incidence scrolls with a
line as directrix curve, those whose base contains a plane,
and finally, those whose base contains a $3$-plane. All the
others may be obtained by  combinations of these
processes.

\subsection{Incidence Scrolls with a Directrix Line}

The following theorem
gives all linearly normal incidence scrolls with a line as
directrix. For the proof we refer the reader to \cite{rosa},
Proposition 4.3.

\begin{teo}\label{p1s}
In $\p^n, \, n\geq 3$, the scroll given by $$\b_n=\{\p^1,
(n-1)\,\p^{n-2}\}$$ in general position is the rational
normal scroll of degree
$n-1$ with a line as minimum directrix.
\end{teo}

We see at once that the previous theorem
gives all rational normal scrolls with a line as
directrix. Then these are projective models of rational ruled
surfaces
$X_e, e\geq 0$. For each $X_e$, the unisecant complete
linear system which gives the immersion is defined by the
very ample divisor 
$H\sim C_o+(e+1)f$.

\begin{center}

\begin{tabular}{|c||c|c|c|c|c|c|c||c||c||c|}
\hline
\hline\multicolumn{11}{|c|}{\small TABLE 1. INCIDENCE SCROLLS
WITH A DIRECTRIX LINE} \\
\hline
\hline
{ }& \multicolumn{7}{|c||}{\scriptsize $n_i, \,
i=1,\cdots ,7$}&{ }&{ }&{ }\\
\cline{2-8} 
{\scriptsize
Scroll}&{\scriptsize
$1$}&{\scriptsize
$2$}&{\scriptsize
$3$}&{\scriptsize
$4$}&{\scriptsize $5$}&{\scriptsize
$6$}&{\scriptsize $7$}&{\scriptsize Normalized}&{\scriptsize
Min. Dir.$({\star})$}&{\scriptsize $deg(\d)$}\\
\hline
{\scriptsize $R^2_0\sub
\p^3$}&{3}&{-}&{-}&{-}&{-}&{-}&{-}&{\scriptsize
$\Te_{\p^1}\oplus\Te_{\p^1}$}&{\scriptsize
$\p^1\;(\infty^1)$ } &{\scriptsize $1$}\\
\hline
\hline
{\scriptsize $R^3_0\sub
\p^4$}&{1}&{3}&{-}&{-}&{-}&{-}&{-}&{\scriptsize
$\Te_{\p^1}\oplus\Te_{\p^1}(-1)$}&{\scriptsize
$\p^1\;(1)$}&{\scriptsize $2$}\\
\hline\hline
{\scriptsize $R^4_0\sub
\p^5$}&{1}&{-}&{4}&{-}&{-}&{-}&{-}&{\scriptsize
$\Te_{\p^1}\oplus\Te_{\p^1}(-2)$}&{\scriptsize
$\p^1\;(1)$}&{\scriptsize $3$}\\
\hline
\hline
{\scriptsize $R^5_0\sub
\p^6$}&{1}&{-}&{-}&{5}&{-}&{-}&{-}&{\scriptsize
$\Te_{\p^1}\oplus\Te_{\p^1}(-3)$}&{\scriptsize
$\p^1\;(1)$}&{\scriptsize $4$}\\
\hline
\hline
{\scriptsize $R^6_0\sub
\p^7$}&{1}&{-}&{-}&{-}&{6}&{-}&{-}&{\scriptsize
$\Te_{\p^1}\oplus\Te_{\p^1}(-4)$}&{\scriptsize
$\p^1\;(1)$}&{\scriptsize $5$}\\
\hline
\hline
{\scriptsize $R^7_0\sub
\p^8$}&{1}&{-}&{-}&{-}&{-}&{7}&{-}&{\scriptsize
$\Te_{\p^1}\oplus\Te_{\p^1}(-5)$}&{\scriptsize
$\p^1\;(1)$}&{\scriptsize $6$}\\
\hline
\hline
{\scriptsize $R^8_0\sub
\p^9$}&{1}&{-}&{-}&{-}&{-}&{-}&{8}&{\scriptsize
$\Te_{\p^1}\oplus\Te_{\p^1}(-6)$}&{\scriptsize
$\p^1\;(1)$}&{\scriptsize $7$}\\
\hline
\hline \multicolumn{11}{l}{$^{\star}$ {\scriptsize Number of
minimum directrix curves}}\\
\end{tabular}

\end{center}

\subsection{Incidence Scrolls with $\p^2 \in \b$}
\label{3}

\begin{teo}\label{p2s}
For every $n\geq 4$ and $0\leq i\leq \frac{n}{2}$, there is
an incidence scroll of degree $d_i^n=(\,
^{n-i}_{\,
\, \, \, 2})+i-1$ and genus $g_i^n=(\, ^{n-i-2}_{\quad 
2})$ given by $$\b_i^n=\{\p^2,
i\,\p^{n-3},(n-2i)\,
\p^{n-2}\}$$ in general position. The directrix curve in
the plane has degree $(n-i-1)$. Moreover, these are all 
incidence scrolls with base $\p^2$.
\end{teo}
{\bf Proof.} If a nondegenerate incidence scroll $R\sub\p^n$ 
has a plane as base space, then the other base spaces have
dimension $n-3$ or $n-2$. Fix $\p^2\in \b$. Then we can vary
the number of 
$\p^{n-3}$'s, written $i$, between $0$ and $\frac{n}{2}$
(respectively
$\frac{n-1}{2}$) if $n$ is even (respectively if $n$ is odd).
Fix $\p^2$ and  $i \, \p^{n-3}$. Then the number of 
$\p^{n-2}$'s is determined by (IS).
 
The proof is by induction on $i$.
If $i=0$, we will proof that the incidence scroll
$R^{d^n_0}_{g^n_0} \sub \p^n$ with base
$\b_0^n=\{\p^2,n\,\p^{n-2}\}$ in general position has degree
$(\, ^{n}_{2})-1$ and genus $(\,^{n-2}_{\,\, \, 2})$, and
the plane directrix curve has degree $(n-1)$. To do
this, we suppose that $\p^2\vee\p^{n-2}=\p^{n-1}\sub\p^n$,
for
$\p^2\in \b_0^n$ and any $\p^{n-2}\in \b^n_0$. Then the
incidence scroll
$R^{d_0^n}_{g_0^n}\sub
\p^n
$ degenerates into:
$R^{d_0^{n-1}}_{g_0^{n-1}}\sub \p^{n-1}$ and the rational
normal scroll of degree $n-1$ with $n-2$ generators in
common. By Proposition \ref{degen}, it is obvious that  
$d_0^n=d_0^{n-1}+n-1$ and $g_0^n=g_0^{n-1}+n-3$. Since the
incidence scroll $R^5_1 \sub\p^4$ is given by $\b_0^4=\{5\,
\p^2\}$ (apply Proposition \ref{degen} to $2\,
\p^2 \in \b^4_0$), 
$$d_0^4=5 \Rightarrow \, d_0^5=9 \Rightarrow \,  d_0^6=14
\Rightarrow \, d_0^n=(\, ^n_2)-1;$$ 
$$g_0^4=1 \Rightarrow \, g_0^5=3 \Rightarrow \, g_0^6=6
\Rightarrow \, g_0^n=(\, ^{n-2}_{\,\,\,2}).$$  
Moreover, the degree
of the plane directrix curve is $(n-1)$ because it is the
number of lines in $\p^n$ which meet a $\p^1$ and $n\,
\p^{n-2}$ (by Giambelli's formula). According to Pieri's
formula, we have:
$$\begin{array}{l}
\om(1, n)\om(n-2, n)^n=\\
=\om(0, n)\om(n-2, n)^{n-1}+\om(1, n-1)\om(n-2, n)^{n-1}=\\
=1+\om(0 ,n-1)\om(n-2, n)^{n-2}+\om(1, n-2)\om(n-2,
n)^{n-2}=\\
 =\cdots=(n-1)\om(0, 1)=n-1.
\end{array}
$$

Assume the theorem holds for $i-1\geq 0$; we will 
obtain the incidence scroll $R^{d_{i-1}^n}_{g_{i-1}^n}\sub
\p^n$ with base
$\b_{i-1}^n=\{\p^2, (i-1)\,\p^{n-3}, (n-2i+2)\, \p^{n-2}\}$
and a directrix curve in $\p^2$ of degree $(n-i)$. 

Separate plane and a
$\p^{n-2}$ in $\b_{i-1}^n$.
Then, we obtain an incidence scroll in $\p^{n+1}$ of
degree
$(\, ^{n-i+1}_{\quad 2})+i-1$ and genus $(\, ^{n-i-1}_{\quad
2})$ with base
$\b '=\{\p^2, i\,\p^{n-2}, (n-2i+1)\, \p^{n-1}\}$, and a 
plane directrix curve of degree $(n-i)$. Hence, writing
$(n-1)$ instead of $n$, we conclude the proof.\qed

\begin{rem} {\em Since we will work under the hypothesis
$n\geq 4$, this theorem gives all the incidence
scrolls with a plane directrix curve. For $n=3$, the
incidence scroll
$R^2_0\sub \p^3$ with a plane directrix curve also appears in
Theorem
\ref{p1s}. Moreover, it is important to note that:
\begin{enumerate}
\item In general, if $R^{d^n_i}_{g^n_i}\sub \p^n \, (n\geq
4)$ has a plane directrix curve, then the plane which
contains such a directrix is a base space (see
\cite{rosa}, Proposition 2.5). This is not true for $n=4$
and $i=2$. In this case, we obtain the incidence scroll
$R^2_0 \sub \p^3$ with base $\b^4_2=\{3\, \p^1\}$ where 
$\p^2 \notin \b^4_2$ because a plane directrix curve is a
hyperplane section and then any plane does not impose
conditions on $G(1, 3)$.
\item The plane directrix curve 
is not necessarily the directrix curve of minimum degree. For
exam\-ple, for
$i=1$ and
$n=4$, we obtain $R^3_0\sub \p^4$ 
which is defined by 
$\b^4_1=\{\p^1, 3\,\p^2\}$ with a
line as minimum directrix curve (see Theorem \ref{p1s}). 
\item For $n=4$, the plane directrix curve is not necessarily
unique. For example, if $i=0$ and
$n=4$, then we obtain $R^5_1\sub \p^4$
with base $\b^4_0=\{5\,\p^2\}$ where each plane directrix
curve is of type 
$C^3_1\sub\p^2$.
\item If we take a plane between the base spaces, then we
impose conditions on the ge\-nus of the incidence
scroll (the same is true for $\p^1$ because the incidence
scroll is always rational). In the notation of Theorem
\ref{p2s} the genus of the scroll (written $g$) is subject
to the condition $$2g=(n-i-2)(n-i-3).$$ For example, for
$g=2$, there exists no incidence scrolls with a plane
directrix curve.
\item There are scrolls with a plane directrix curve 
which are not incidence. For example, we can take $R^6_0\sub
\p^7$ with a directrix conic, a three-dimensional family of
directrix quartics and a five-di\-men\-sional family of
directrix quintics. It is not defined by incidences (see
\cite{rosa}, Example 4.8).
\end{enumerate}
}
\end {rem}

\begin{center}

\begin{tabular}{|c||c|c|c|c|c|c|c||c|}
\hline
\hline\multicolumn{9}{|c|}{\scriptsize TABLE 2. INCIDENCE
SCROLLS WITH A BASE PLANE} \\
\hline
\hline
{ }& \multicolumn{7}{|c||}{\scriptsize $n_i, \,
i=1,\cdots ,7$}&{ }\\
\cline{2-8} 
{\scriptsize
Scroll}&{\scriptsize
$1$}&{\scriptsize
$2$}&{\scriptsize
$3$}&{\scriptsize $4$}&{\scriptsize
$5$}&{\scriptsize $6$}&{\scriptsize $7$}&{\scriptsize Min.
Dir.$(\star
\star)$}\\
\hline
\hline\multicolumn{9}{||l||}{\bf Genus 0 $\;
(n=i+3) \Rightarrow (n \leq 6)$}\\
\hline
\hline
{\scriptsize $R^2_0\sub
\p^3$}&{3}&{-}&{-}&{-}&{-}&{-}&{-}&{\scriptsize
$\p^1\;(\infty^1)$}\\
\hline
\hline
{\scriptsize $R^3_0\sub
\p^4$}&{1}&{3}&{-}&{-}&{-}&{-}&{-}&{\scriptsize
$\p^1\;(1)$}\\
\hline\hline
{\scriptsize $R^4_0\sub
\p^5$}&{-}&{3}&{1}&{-}&{-}&{-}&{-}&{\scriptsize $C^2_0\sub
\p^2\;(\infty^1)$}\\
\hline
\hline
{\scriptsize $R^5_0\sub
\p^6$}&{-}&{1}&{3}&{-}&{-}&{-}&{-}&{\scriptsize $C^2_0\sub
\p^2\;(1)$}\\
\hline
\hline\multicolumn{9}{||l||}{\bf Genus 1 $\;
(n=i+4)\Rightarrow (n \leq 8)$}\\
\hline
\hline
{\scriptsize $R^5_1\sub
\p^4$}&{-}&{5}&{-}&{-}&{-}&{-}&{-}&{\scriptsize $C^3_1\sub
\p^2\;(\infty^1)$}\\
\hline\hline
{\scriptsize $R^6_1\sub
\p^5$}&{-}&{2}&{3}&{-}&{-}&{-}&{-}&{\scriptsize $C^3_1\sub
\p^2\;(2)$}\\
\hline
\hline
{\scriptsize $R^7_1\sub
\p^6$}&{-}&{1}&{2}&{2}&{-}&{-}&{-}&{\scriptsize $C^3_1\sub
\p^2\;(1)$}\\
\hline
\hline
{\scriptsize $R^8_1\sub
\p^7$}&{-}&{1}&{-}&{3}&{1}&{-}&{-}&{\scriptsize $C^3_1\sub
\p^2\;(1)$}\\
\hline
\hline
{\scriptsize $R^9_1\sub
\p^8$}&{-}&{1}&{-}&{-}&{4}&{-}&{-}&{\scriptsize $C^3_1\sub
\p^2\;(1)$}\\
\hline
\hline\multicolumn{9}{||l||}{\bf Genus 3 $\;
(n=i+5)\Rightarrow (n \leq 10)$}\\
\hline
\hline
{\scriptsize $R^9_3\sub
\p^5 \star$}&{-}&{1}&{5}&{-}&{-}&{-}&{-}&{\scriptsize
$C^4_3\sub \p^2\;(1)$}\\
\hline\hline
{\scriptsize $R^{10}_3\sub
\p^6\star$}&{-}&{1}&{1}&{4}&{-}&{-}&{-}&{\scriptsize
$C^4_3\sub \p^2\;(1)$}\\
\hline
\hline
{\scriptsize $R^{11}_3\sub
\p^7\star$}&{-}&{1}&{-}&{2}&{3}&{-}&{-}&{\scriptsize
$C^4_3\sub \p^2\;(1)$}\\
\hline
\hline
{\scriptsize $R^{12}_3\sub
\p^8\star$}&{-}&{1}&{-}&{-}&{3}&{2}&{-}&{\scriptsize
$C^4_3\sub \p^2\;(1)$}\\
\hline\hline
{\scriptsize $R^{13}_3\sub
\p^9\star$}&{-}&{1}&{-}&{-}&{-}&{4}&{1}&{\scriptsize
$C^4_3\sub \p^2\;(1)$}\\
\hline
\hline
{\scriptsize $R^{14}_3\sub
\p^{10}\star$}&{-}&{1}&{-}&{-}&{-}&{-}&{5}&{\scriptsize
$C^4_3\sub \p^2\;(1)$}\\
\hline
\hline\multicolumn{9}{l}{{\scriptsize ${\star}$} {
Incidence special scroll}}\\ \multicolumn{9}{l}{{\scriptsize
$(\star
\star)$} {Number of minimum directrix curves}}\\
\end{tabular}
\end{center}

In Table 2 we see some examples of
incidence scrolls which have a plane directrix curve.
Moreover, the table contains all the incidence scrolls with a
plane directrix curve of genus $0, 1$ and
$3$.

\subsection{Incidence Scrolls with $\p^r \in \b, r\geq 3$}
\label{4}

We can continue with a similar method for obtain a
new theorem which provides all the incidence scrolls which
have a $\p^3$ between the base spaces. Without loss of
generality we can assume $n \geq 5$ because the other cases
appear in Theorems \ref{p1s} and \ref{p2s}. For example,
$R^3_0\sub\p^4$ with $\p^1$ as minimum directrix curve is
given by Theorem \ref{p1s}.

Fix $\p^3$. The other base spaces have dimension $\geq n-4$
because the scroll is non-degenerate. Then the other base
spaces are necessarily
$j\,\p^{n-4}$'s,
$i\,
\p^{n-3}$'s and
$n+1-3j-2i\, \p^{n-2}$'s. In fact, we see at once
that
\begin{enumerate}
\item $\p^3$ imposes $(n-4)$ conditions but we need $(2n-3)$
conditions;
\item $0\leq j\leq \frac{1}{3}(n+1)$;
\item fixing j between the above limits, we obtain  
$0\leq i \leq \frac{1}{2}(n+1-3j)$.
\end{enumerate}

\begin{lemma}\label{p3s1} 
In $\p^n \, (n\geq 5)$, the incidence scroll given by 
$\b^n_{0,0}=\{\p^3, (n+1)\, \p^{n-2}\}$ in general
position  has degree
$d^n_{0,0}=(\, ^{n+1}_{\,
\, \, 3})-n-1$ and genus $g^n_{0,0}=(\,
^{n}_{3})+(\,^{n-1}_{\, \, \, 3})-2n+4$. The directrix
curve in $\p^3$ has degree
$(\, ^{n}_{2})-1$.
\end{lemma}
{\bf Proof.} We proceed by induction on $n$.
If $n=5$, then we will prove that $\b^5_{0,0}$
defines the scroll $R^{14}_8\sub \p^5$ and the
directrix curve in $\p^3$ has degree $9$. Suppose that $2\,
\p^3$'s have by intersection a plane instead of a line.
Since the scroll breaks up into $R^9_3\sub 
\p^5$ and $R^5_1\sub \p^4$ (see Theorem \ref{p2s}) with $5$
generators in common, Proposition \ref{degen} shows that the
scroll has degree $14$ and genus $8$. Moreover, the
directrix curve in each $\p^3$ has degree equal to the number
of lines in $\p^5$ which meet a $\p^2$ and $6\,\p^3$'s. From
Pieri's formula we have:
$$\begin{array}{l}
\om(2, 5)\om(3, 5)^6=\\
=\om(1, 5)\om(3, 5)^{5}+\om(2, 4)\om(3,5)^{5}=\\ 
=\om(0, 5)\om(3, 5)^4+2\om(1, 4)\om(3, 5)^{4}+\om(2, 3)\om(3,
5)^{4}=\\
=\cdots=9\om(0, 2)\om(3, 5)=9\om(0, 1)=9.
\end{array}
$$

Suppose that the incidence scroll
with base $\b^n_{0,0}$ in general position has degree
$d^n_{0,0}$ and genus $g^n_{0,0}$ and each directrix
curve in $\p^3$ has degree $(\, ^{n}_{2})-1 $. Then
we will proof the lemma for $n+1$. 

Let $R^{d}_{g}\sub \p^{n+1}$ be the incidence
scroll defined by $\b^{n+1}_{0,0}$ in general position. We
will calculate its degree and genus. If
$\p^3\vee\p^{n-1}=\p^{n}$, then the scroll breaks up into:

\begin{enumerate}
\item[-] $R^{(\,^{n+1}_{\,\, \, 2})-1}_{(\,^{n-1}_{\,\, \,
2})}\sub \p^{n+1}$ with base $\b^{n+1}_{0}$ (see Theorem
\ref{p2s});
\item[-] $R^{d^n_{0,0}}_{g^n_{0,0}}\sub \p^n $ with base
$\b^{n}_{0,0}$, by hypothesis of induction, and $(\,
^{n+1}_{\, \,
\, 2})-n-1$ generators in common.
\end{enumerate}

An easy computation shows that $d=(\, ^{n+2}_{\, \, \,
3})-n-2=d^{n+1}_{0,0}$ and $g=(\,
^{n}_{3})+(\,^{n+1}_{\,
\, \, 3})-2n+2=g^{n+1}_{0,0}$. Moreover, the degree of the
directrix curve in $\p^3$ is the number of lines in
$\p^{n+1}$ which meet a $\p^2$ and 
$(n+2)\, \p^{n-1}$'s. It is exactly:
$$\begin{array}{l}
\om(2, n+1)\om(n-1, n+1)^{n+2}=\\
=\om(1, n+1)\om(n-1, n+1)^{n+1}+\om(2, n)\om(n-1,
n+1)^{n+1}=\\  
=n\om(0, 1)+(n-1)\om(0, 1)+\om(2, n-1)\om(n-1, n+1)^n=\\
=\cdots=n+(n-1)+(n-2)+\cdots+2=\frac{1}{2}(n-1)(n+2).\qed
\end{array}
$$

\begin{lemma}\label{p3s2} 
In $\p^n \, (n\geq 5)$, the incidence scroll with base
$\b^n_{0,1}=\{\p^3, \p^{n-3},(n-1)\,\p^{n-2}\}$ in general
position has degree
$d^n_{0,1}=(\, ^{n}_{3})-1$ and genus $g^n_{0,1}=(\,
^{n-2}_{\, \, 3})+(\,^{n-1}_{\,
\, 3})-n+3$. The directrix curve in $\p^3$ has degree
$(^{n-1}_{\, \, \, \, 2})$. 
\end{lemma}
{\bf Proof.}
Suppose that $\p^3\vee\p^{n-3}=\p^{n-1}$. Then the incidence 
scroll $R^d_g\sub \p^n$ with base $\b^n_{0,1}$ degenerates
into:

\begin{enumerate}
\item[-] $R^{d^{n-1}_{0,0}}_{g^{n-1}_{0,0}}\sub\p^{n-1}$
with base $\b^{n-1}_{0,0}$ (see Lemma
\ref{p3s1});
\item[-] $R^{n-1}_{0}\sub 	\p^n$ with base $\b_{n-1}=\{\p^1,
(n-1)\,\p^{n-2}\}$ (see Theorem \ref{p1s}) and $n-2$
generators in  common.
\end{enumerate}

Whence, we conclude that $d=(\,^{n}_{3})-1=d^n_{0,1}$
and $g=(\, ^{n-2}_{\,\,\, 3})+(\,^{n-1}_{\, \, \,
3})-n+3=g^n_{0,1}$.

The degree of the directrix curve which is contained in
$\p^3$ is given by:
$$\begin{array}{l}
\om(2, n)\om(n-3, n)\om(n-2, n)^{n-1}=\\
=\om(0, n)\om(n-2, n)^{n-1}+\om(1, n-1)\om(n-2,
n)^{n-1}+\\+\om(2, n-2)\om(n-2, n)^{n-1}=\cdots=\\
=1+(n-2)+(n-3)+\cdots+2=\frac{1}{2}(n-1)(n-2).
\end{array}
$$\qed

\begin{lemma} \label{p3s3}
In $\p^n \, (n\geq 5)$, the incidence scroll with base
$\b^n_{1,0}=\{\p^3, \p^{n-4}, (n-2)\, \p^{n-2}\}$ in
general position has degree
$d^n_{1,0}=(\, ^{n-1}_{\, \, \, 3})$ and genus $g^n_{1,0}=(\,
^{n-2}_{\, \, 3})+(\,^{n-3}_{\, \, 3})-n+4$. The directrix
curve in $\p^3$ has degree
$(^{n-2}_{\, \, \, \, 2})$.
\end{lemma}
{\bf Proof.}
Separate $\p^3$ and $\p^{n-4}$. Then we conclude from
Proposition \ref{degen} that our scroll
$R^d_g\sub \p^n$ with base $\b^n_{1,0}$ degenerates into
$R^{d-1}_{g}\sub
\p^{n-1}$ with base $\b^{n-1}_{0,1}$, hence that
$d=(\,^{n-1}_{\,\,\, 3})$ and $g=(\, ^{n-2}_{\, \,
3})+(\,^{n-3}_{\, \, 3})-n+4$, by Lemma
\ref{p3s2}.\qed
\bigskip

If we work these lemmas, then we can obtain a
general theorem which will give every incidence scroll
with a $\p^3$ between the base spaces. Moreover, we
will obtain the degree and the genus of the incidence scroll
for any
$i$ and $j$ between suitable limits. The degree of the
directrix curve of such an incidence scroll which is given
by the following theorem may easily be found by referring to
each particular case and using Pieri's formula.

\begin{teo} \label{p3s}
For every $n\geq 5, \, 0\leq j\leq \frac{1}{3}(n+1)$
and $0\leq i \leq \frac{1}{2}(n+1-3j)$, there is an incidence
scroll of degree
$$d^n_{j,i}=d^{n-j}_{0,i+j}+j=(^{n-i-2j+1}_{\quad \, \, \,
3})-(n-i-2j)+(i+j)(n-i-2j-1)+j-1$$ and genus
$$g^n_{j,i}=g^{n-j}_{0,i+j}=(^{n-i-2j}_{\quad
3})+(\,^{n-i-2j-1}_{\quad \, \,
\, 3})-2(n-i-2j)+(i+j)(n-i-2j-2)+4$$ with base
$\b^n_{j,i}=\{\p^3, j\, \p^{n-4},i\,\p^{n-3},
(n+1-3j-2i)\,\p^{n-2}\}$ in general position.
\end{teo}
{\bf Proof.} We proceed by induction on $j$.

For $j=0$, we will proof that the incidence
scroll with base $\b^n_{0,i}$, in general position, has
degree $d^n_{0,i}=(^{n-i+1}_{\quad 
3})-(n-i)+i(n-i)-1$ and genus
$g^n_{0,i}=(^{n-i}_{\, \, \,
3})+(\,^{n-i-1}_{\quad 3})-2(n-i)+i(n-i-2)+4$.
To do this, we proceed by induction on $i$. The theorem is
true for $i=0$, by Lemma
\ref{p3s1}. Supposing the theorem true for $i-1\geq 0$, we
prove it for
$i$. Suppose that $\p^3\vee\p^{n-3}=\p^{1}$, for any
$\p^{n-3}\in \b^n_{0,i}$ . Then the incidence  scroll
$R^d_g\sub \p^n$ with base $\b^n_{0,i}$ in general position
degenerates
into: $R^{d^{n-1}_{0,i-1}}_{g^{n-1}_{0,i-1}}\sub\p^{n-1}$
with base $\b^{n-1}_{0,i-1}$( by hypothesis of induction) and
$R^{n-i}_{0}\sub 	\p^{n-i+1}$ with base
$\b_{n-i}=\{\p^1, (n-i)\,\p^{n-i-1}\}$ (see Theorem
\ref{p1s}) and $(n-i-1)$ generators in  common.

Whence, we find that $d=d^{n-1}_{0,i-1}+(n-i)=d^n_{0,i}$
and $g=g^{n-1}_{0,i-1}+(n-i-2)=g^n_{0,i}$.

Supposing the theorem true for $j-1\geq 0$, we prove it for
$j$. Separating $\p^3$ and $\p^{n-3}$ in $\b^n_{j-1,i}$ we
obtain a base $\b=\{\p^3, j\, \p^{n-3},$ $ (i-1)\,$ $
\p^{n-2}, (n+1-3(j-1)-2i)\,
\p^{n-1}\}$ in general position which defines
$R^{d^{n}_{j-1,i}+1}_{g^{n}_{j-1,i}}\sub\p^{n+1}$. 

Write
$n$ instead $n+1$ and $i$ instead $i-1$. Then the incidence
scroll $R^d_g \sub \p^n$ with base $\b^n_{j,i}$ in general
position has degree $d=d^{n-1}_{j-1,i+1}+1$ and genus
$g=g^{n-1}_{j-1,i+1}$. Therefore, $d=d^{n-j}_{0,i+j}+j$
and $g=g^{n-j}_{0,i+j}$, by hypothesis of induction.\qed

\begin{rem}
{\em The determination of the degree of the directrix curve
in $\p^3$, written $deg(C)$, is more subtle than in the
other cases because it is $\om(2, n)\om(n-4, n)^j\om(n-3,
n)^i\om(n-2, n)^{n+1-3j-2i}$. From what has already been
proved it may be concluded that for every $j,i$ between
suitable limits,
$$deg(C)=(^{n-i-2j} _{\, \, \, \, \, \, \, 2})+i+j-1.$$
}
\end{rem}

If we now apply this argument again, with $\p^3$ replaced by
$\p^r, \, r\geq4$, then we can obtain a general theorem which
will give all incidence scrolls with $\p^r$ between
the base spaces. Its formulation is very complicated because
we must work with, at least, three index: $i,j,k$. But
the proof is similar to that of Theorem \ref{p3s} for all
cases. Therefore we can obtain degree and genus of the
incidence scroll with base
$$\begin{array}{r}
\b^n_{i_1,\cdots,i_{r-1}}=\{\p^r, i_{1}\,
\p^{n-r-1},\cdots, i_{r-2}\,
\p^{n-4}, i_{r-1}\,\p^{n-3},\\
(n+r-2-ri_1-\cdots-3i_{r-2}-2i_{r-1})\,
\p^{n-2}\}
\end{array}$$ 
for any $i_1, \cdots , i_{r-1}$ between
suitable limits.

\begin{center}

\begin{tabular}{|c||c|c|c|c||c|}
\hline
\hline\multicolumn{6}{||c||}{\scriptsize TABLE 3. INCIDENCE
SCROLLS WITH $\p^3\in \b$} \\
\hline
\hline
{ }& \multicolumn{4}{|c||}{\scriptsize $n_i, \,
i=1,\cdots ,4$}&{ }\\
\cline{2-5} 
{\scriptsize
Scroll}&{\scriptsize
$2$}&{\scriptsize $3$}&{\scriptsize
$4$}&{\scriptsize $5$}&{\scriptsize Directrix in $\p^3$}\\
\hline
\hline
{\scriptsize $R^{14}_8\sub
\p^5\star$}&{-}&{7}&{-}&{-}&{\scriptsize $C^9_8\sub \p^3$}\\
\hline\hline
{\scriptsize $R^5_0\sub
\p^6$}&{1}&{3}&{-}&{-}&{\scriptsize $C^3_0\sub \p^3$}\\
\hline
{\scriptsize $R^7_1\sub
\p^6$}&{1}&{2}&{2}&{-}&{\scriptsize $C^4_1\sub \p^3$}\\
\hline
{\scriptsize $R^{10}_3\sub
\p^6\star$}&{1}&{1}&{4}&{-}&{\scriptsize $C^5_3\sub \p^3$}\\
\hline
{\scriptsize $R^9_2\sub
\p^6$}&{-}&{4}&{1}&{-}&{\scriptsize $C^5_2\sub \p^3$}\\
\hline
{\scriptsize $R^{13}_5\sub
\p^6 \star$}&{}&{3}&{3}&{-}&{\scriptsize $C^7_5\sub \p^3$}\\
\hline
{\scriptsize $R^{19}_{11}\sub
\p^6 \star$}&{-}&{2}&{5}&{-}&{\scriptsize $C^{10}_{11}\sub
\p^3$}\\
\hline
{\scriptsize $R^{28}_{22}\sub
\p^6\star$}&{-}&{1}&{7}&{-}&{\scriptsize $C^{14}_{22}\sub
\p^3$}\\
\hline
\hline
{\scriptsize $R^6_0\sub
\p^7$}&{-}&{3}&{1}&{-}&{\scriptsize $C^3_0\sub \p^3$}\\
\hline
{\scriptsize $R^{8}_1\sub
\p^7$}&{-}&{3}&{-}&{2}&{\scriptsize $C^4_1\sub \p^3$}\\
\hline
{\scriptsize $R^{10}_2\sub
\p^7$}&{-}&{2}&{2}&{1}&{\scriptsize $C^5_2\sub \p^3$}\\
\hline
{\scriptsize $R^{14}_5\sub
\p^7\star$}&{-}&{2}&{1}&{3}&{\scriptsize $C^7_5\sub \p^3$}\\
\hline
{\scriptsize $R^{20}_{11}\sub
\p^7\star$}&{-}&{2}&{-}&{5}&{\scriptsize $C^{10}_{11}\sub
\p^3$}\\
\hline
{\scriptsize $R^{12}_3\sub
\p^{7}$}&{-}&{1}&{4}&{-}&{\scriptsize $C^6_3\sub \p^3$}\\
\hline
\hline\multicolumn{6}{l}{$^{\star}$ {\scriptsize Incidence
special scroll}}
\end{tabular}
\end{center}


\begin{thebibliography}{77}

\addcontentsline{toc}{section}{Bibliography}

\bibitem{rosa}{\sc  CID-MU$\tilde N$OZ,R.-PEDREIRA,M.}
{\it Classification of Incidence Scrolls(I); preprint
math.AG/0005272.}\par  
{To appear in Manuscripta Mathematica}.


\bibitem{gyh}{\sc  GRIFFITHS, P. - HARRIS, J. }
{\it Principles of algebraic geometry.}\par { Pure and
applied mathematics. }  A Wiley-Interscience publication,
pp. 193-211 (1978).


\bibitem{hatshor}{\sc  HARTSHORNE, R. }
{\it Algebraic Geometry.}\par { Graduate Texts in 
Mathematics 52. }  Springer-Verlag, New York. 1977.


\bibitem{Kleiman}{\sc KLEIMAN, S. L. - LAKSOV, D. }
{\it Schubert calculus.}\par
{ Amer. Math. Monthly 79, }  pp. 1061-1082 (1972).


\end{thebibliography}
\end{document}